\documentclass[12pt]{amsart}
\usepackage{amsmath, amsfonts, amssymb, amscd}

\newtheorem{theorem}{Theorem}[section]

\newtheorem{lemma}[theorem]{Lemma}
\newtheorem{corollary}[theorem]{Corollary}

\theoremstyle{definition}

\newcommand\R{{\mathbb R}}
\newcommand\C{{\mathbb C}}
\newcommand\aut{\mathrm{Aut}}
\newcommand\GL{\mathrm{GL}}

\newcommand\Sp{\mathrm{Sp}}
\newcommand\Or{\mathrm{O}}
\newcommand\alt{\mathrm{Alt}}
\newcommand\sym{\mathrm{Sym}}

\title[Permutation Polytopes]{Permutation Polytopes and Indecomposable 
Elements in
Permutation Groups}
\author{Robert Guralnick and David Perkinson}
\date{\today}

\address
{Robert M. Guralnick\\
Department of Mathematics\\
University of Southern California\\
3620 S. Vermont Ave.\\
Los Angeles CA 90089-2532}
\email{\tt guralnic@usc.edu}

\address{David Perkinson\\
Department of Mathematics\\
Reed College\\
Portland, OR 97202}
\email{\tt davidp@reed.edu}

\thanks{The first author acknowledges support from
the NSF grant DMS 0140578. We thank D. Goldstein
for his comments on earlier drafts. We also thank J. Buhler
and I. Pak for their input. Finally, we thank the referees
for their comments.}

\begin{document}

\begin{abstract}
Each group $G$ of $n\times n$ permutation matrices has a corresponding
permutation polytope, $P(G):=\mbox{conv}(G)\subset\R^{n\times n}$.  We
relate the structure of $P(G)$ to the transitivity of $G$.  In particular, we
show that if $G$ has $t$ nontrivial orbits, then $\min\{2t,\lfloor
n/2\rfloor\}$ is a sharp upper bound on the diameter of the graph of
$P(G)$.  We also show that $P(G)$ achieves its maximal dimension of
$(n-1)^2$ precisely when $G$ is $2$-transitive.  We then extend the
results of Pak \cite{Pa} on mixing times for a random walk on $P(G)$.
Our work depends on a new result for permutation groups involving
writing permutations as products of indecomposable permutations.
\end{abstract}

\maketitle

\section{Introduction}
Let $G$ be a subgroup of $S_n$, the symmetric group on
$\{1,2,\dots,n\}$. Via the usual representation of $G$ as a group of
$n\times n$ permutation matrices, each element of $G$ may be
considered as an element of $\R^{n^2}$. The convex hull in $\R^{n^2}$
of the elements of $G$ is $P(G)$, the {\em permutation polytope}
associated with $G$. Permutation polytopes and their linear
projections have been studied extensively due to their connection to
problems in combinatorial optimization \cite{Ba}, \cite{BV},
\cite{On}, \cite{YKK}. The most well-known example is the case where
$G=S_n$ with corresponding permutation polytope called the $n$-th {\em
Birkhoff polytope} or the $n$-th {\em assignment polytope} \cite{Bi},
\cite{BG}, \cite{YKK}. Even here there are open problems
\cite{Pa}; for instance, its volume is known only up to $n=10$
\cite{BP}. Some newer applications of permutation polytopes are to
group resolutions \cite{El} and communications networks \cite{LJD},
\cite{Ti}.

The main concern of this paper is to establish links between
(algebraic) properties of an arbitrary permutation group
$G$ and (geometric) properties of its corresponding permutation
polytope $P(G)$. We are especially interested in ways in which
the transitivity of $G$ is reflected in its polytope.
First, Theorem~\ref{weaker}, shows
that every element of a transitive permutation group can be written as
a product of at most two so-called indecomposable elements (see
\S\ref{Permutation Groups} for definitions). The geometric
consequence is Corollary~\ref{diameter}: if $G$ has $t$ non-trivial
orbits, then the diameter of $P(G)$, i.e., the diameter of the edge
graph of $P(G)$, is bounded by $\min\{2t,\lfloor n/2\rfloor\}$.
Thus, if $G$ is transitive, the diameter of $P(G)$ is at most $2$.
This generalizes previous work establishing the diameters of the
Birkhoff polytopes \cite{BR}, \cite{Yo} and the diameters of the
polytopes corresponding to the groups of even permutations
\cite{BL}. In the language of Babai et al. \cite{BHKLS}, we have
bounded the diameter of {\em the group} $G$ with respect to the set of
generators consisting of its indecomposable elements.

Corollary~\ref{diameter} relies on Theorem~\ref{thm:smallest face}
characterizing the smallest face of a permutation polytope
containing two prescribed vertices (group elements) in terms of
their cycle structure. In particular, we characterize the edges of
a permutation polytope, as previously known for the Birkhoff
polytopes \cite{PR} and for the polytopes corresponding to the
groups of even permutations \cite{BL}. The special case $G=S_n$
in Theorem~\ref{thm:smallest face} is Proposition~2.1 in \cite{BS}.

The other main result concerning transitivity is
Corollary~\ref{cor:transitivity}, showing that the
dimension of $P(G)$ is bounded by $(n-1)^2$ with equality if and only if $G$ is
$2$-transitive. The dimension of the $n$-th Birkhoff polytope is
known to equal the maximum value, $(n-1)^2$, by an easy calculation in
linear algebra. With more work, one may similarly show that the
maximum dimension is achieved when $G$ is the collection of all even
permutations and $n\geq 4$ \cite{BL}.
Corollary~\ref{cor:transitivity} generalizes these results and
provides a conceptual explanation.

In the final section of the paper, we generalize the results of
\cite{Pa} about the mixing time of random walks on these
polytopes. This says that random products of indecomposable elements
tend to the uniform distribution very quickly for $G$ primitive (Pak
\cite{Pa} handles the case of the Birkhoff polytope).

The results in this paper stem from systematic experimentation using
the computer programs GAP \cite{Ga} for group theory and Polymake
\cite{Po} for polytopes.

\section{Permutation Groups}\label{Permutation Groups}
Let $G$ be a permutation group acting faithfully on a (finite) set
$X$. We say $g \in G$ is indecomposable if $g \ne xy$ where $x,y$ are
nontrivial elements of $G$ and $M(x) \cap M(y)$ is empty, where $M(x)$
is the {\em support} of $x$: the set of points of $X$ moved by
$x$. Let $F(x)$ be the set of fixed points of $x$ and $f(x)=|F(x)|$.

We shall prove:

\begin{theorem}\label{weaker} Let $G$ be transitive
on $X$. Then every element of $G$
is a product of at most $2$ indecomposable
elements.
\end{theorem}

In fact, for inductive purposes, it is better to prove
a slightly stronger result:

\begin{theorem} \label{stronger}
Let $G$ be transitive
on $X$. Then every element of $G$
is a product of two elements, each indecomposable
and at least one fixed point free.
\end{theorem}

We will prove this result in the next few subsections.
We first show that it suffices to assume that $G$ acts
primitively on the set $X$ (i.e. preserves no nontrivial
partition of $X$).

We then show that the result holds when the group is
primitive and not almost simple (recall a group
is almost simple if it has a unique minimal normal
subgroup that is a nonabelian simple group).

Finally, we show that in the almost simple
case, aside from the case that $G$ contains $\alt(X)$,
every element is indecomposable (whence the result
follows since fixed point free elements in a finite
transitive permutation group always exist).
The result in the case $G=\alt(X)$ or $\sym(X)$ is
elementary.

We do have to invoke the classification of finite
simple groups to handle the case that $G$ is almost
simple. The key result we use is the classification
of primitive permutation groups containing a nontrivial element
with $f(x) \ge |X|/2$.

We first point out some easy consequences of Theorem~\ref{stronger}
using the following lemma.

\begin{lemma} Suppose that $X = Y \cup Z$ is a finite
$G$-set with $Y$ and $Z$ invariant under $G$. Let
$N$ be the normal subgroup of $G$ acting trivially on
$Y$. If every element of $G/N$ acting on
$Y$ can be written as a product of
$r$ indecomposables and every element of $N$ can
be written as a product of $s$ indecomposables,
then every element of $G$ is a product of
$r+s$ indecomposables.
\end{lemma}

\begin{proof} If $g \in G$, let $g_Y$ denote $g$
considered as permutation on $Y$.

We claim that if $g \in G$ and $g_Y$ is indecomposable,
then $gn$ is indecomposable for some $n \in N$.

Proof of Claim: If $g$ is indecomposable, we are done.
If not, write $g=hu$ where $M(h) \cap M(u)$ is empty
and $h$ is not in $N$. Since $g_Y$ is indecomposable,
$h_Y=g_Y$ and $u \in N$. Thus, $h \in gN$ is indecomposable.

The claim implies that we can write $g \in G$
as a product of $r$ indecomposables (or fewer)
times an element of $N$. By assumption, the
element in $N$ can be written as a product of
$s$ indecomposables (in $N$ and thus also in $G$).
\end{proof}

\begin{corollary}\label{cor:bound1}
If $G \le S_n$, then every element of $G$
can be written as a product of $2t$ indecomposables
where $t$ is the number of nontrivial orbits of $G$.
\end{corollary}

\begin{corollary}\label{cor:bound2} If $G \le S_n$,
then every element of $G$
can be written as a product of $\lfloor n/2\rfloor$ indecomposables.
\end{corollary}

\begin{proof} By induction and the lemma above, it suffices
to consider the case that $G$ is transitive. By the theorem,
the result holds for $n \ge 4$. Inspection shows that
for $n \le 3$, every nontrivial element is indecomposable.
\end{proof}

\subsection{Reduction to the Primitive Case}

Let $G$ be a group acting faithfully and transitively
on the finite set $X$. Let $n=|X| > 1$.

\begin{lemma} Let $Y:=\{X_1, \ldots, X_m\}$ be a nontrivial
$G$-invariant
partition of $X$. Let $N$ be the normal subgroup of $G$
preserving each $X_i$.  Let $g\in G$.
\begin{enumerate}
\item If $gN$ is fixed point free and indecomposable on $Y$,
then every element in $gN$ is fixed point free and indecomposable on
$X$.
\item If $gN$ is indecomposable on $Y$, then there is some element
in $gN$ that is indecomposable on $X$.
\end{enumerate}
\end{lemma}

\begin{proof} We prove both statements simultaneously.
Reordering if necessary, we may assume that $g$ moves
the sets $X_1, \ldots, X_e$ and fixes the other $X_i$.
Assume also that $gN$ is indecomposable on $Y$.

Suppose that $g=xy$ where $M(x) \cap M(y)$ is empty.
Then $gN=xNyN$ and $xN$ and $yN$ cannot move a common $X_i$.
Since $gN$ is indecomposable, we may assume that $gN=xN$
and $yN=N$. Thus, $x \in gN$ and the second statement holds.

Moreover, since $x$ and $y$ share no moved points, $y$
must be trivial on each block moved by $g$. So if $gN$ has no fixed
points on $Y$, then $y=1$ and $g=x$ is indecomposable.
\end{proof}

An immediate consequence is:

\begin{corollary} \label{reductiontoprimitive}
Suppose that $(G,X)$ is a counterexample to Theorem \ref{stronger}
with $|X|$ minimal.  Then $G$ acts primitively on $X$.
\end{corollary}

\begin{proof}  If $G$ preserves a nontrivial partition $Y$
on $X$, let $N$ be the normal subgroup acting trivially
on the partition.  By the previous result, $(G/N,Y)$
is a counterexample to Theorem \ref{stronger}, contradicting
the minimality of $|X|$.
\end{proof}
 
We deal with the case  that $G$ acts primitively on 
$X$ in the next two subsections.

\subsection{Primitive Groups I}

In this subsection, we assume that $G$ is not almost simple
and acts primitively (and faithfully) on the finite set
$X$ of cardinality $n$.

The structure of finite primitive groups is quite constrained.
See \cite{AS} for a detailed description.

Recall that a transvection is a nontrivial unipotent linear
transformation which is trivial on a hyperplane.

\begin{theorem} \label{regnormal} Assume that $G$ contains
a regular normal subgroup $N$. Then one of the following
holds:
\begin{enumerate}
\item Every element of $G$ is indecomposable.
\item $N$ is an elementary abelian $2$-group of order
$2^a \ge 4$ and $G=NH$
where $H$ is a subgroup of $\GL(a,2)=\aut(N)$ acting
irreducibly on $N$ and containing transvections.
\end{enumerate}
Moreover, $G$ satisfies the conclusion of
Theorem \ref{stronger}.
\end{theorem}

\begin{proof} It follows by \cite{AS} that $N$
is a direct product of isomorphic copies of a simple
group $L$. If $g \in G$ has a fixed point, then
as $g$-set, we can identify $X$ with $N$ and the
fixed points of $g$ are identified with $C_N(g)$.
Unless $|L|=2$, any proper subgroup of $N$ has index
at least $3$,
so for $1 \ne g$, the proportion of fixed points is
at most $1/3$. Thus, $M(x) \cap M(y)$ is nonempty for
any two nontrivial elements in $G$ and so (1) holds.

So $N$ is an elementary abelian $2$-group of order $2^a$.
If $a=1$, then $G$ is cyclic of order $2$ and the result
hold. If $a > 1$, the argument of the previous paragraph
applies and we see that $f(g) \le n/2$ with equality if
and only if $g$ induces a transvection acting on $N$.
Thus either (1) or (2) hold.

So it suffices to prove the last statement in the
case $G=NH$ where $|N|=n=2^a \ge 4$ and $G=NH$ with
$H$ acting irreducibly and faithfully on $N$ and containing
transvections. Note that if $x \in G$ is decomposable,
then $x=uv$ where $u,v$ are involutions fixing precisely
one half the points of $X$. Moreover, $u$ and $v$
commute and the fixed point sets of $u$ and $v$ must be
disjoint. Thus, $x$ is a fixed point free involution.

If $a=2$, then $G=S_4$ and the result holds by inspection.
So assume that $a > 2$ and $g$ is a fixed point free involution.

First suppose that $g \in N$.  Choose $h_1, h_2 \in H$ that
are noncommuting transvections (if all transvections in $H$ commute
they would generate a normal unipotent subgroup of $H$ and this
contradicts the irreducibility of $H$).   So $h_1h_2$ has order 
$3$ and  $\langle h_1, h_2 \rangle$ 
centralizes a subgroup $N_0$ a  subgroup of index $4$ in $N$.
Let $1 \ne v \in N_0$ (this is possible since $a > 2$).  Then
$h:=h_1h_2v$ has order $6$ and is fixed point free (since $h^3=v$ is).
Finally, we see that $g=h(h^{-1}g)$ and $h^{-1}g$ has order a multiple
of $3$ and so is indecomposable.

Finally, suppose that $g$ is a fixed point free involution not in $N$.
Let $h_1$ and $h_2$ be noncommuting transvections in $H$.   
Choose $v_i \in N, 1 \le i \le 2$ so that $w_i:=h_iv_i$ has order $4$
(and so is fixed point free and indecomposable).  Let $v$ be a nontrivial
element of $N_0$ (as in the previous paragraph).  Set $w_3:=h_1h_2v$.
So $w_3$ has order $6$ and is fixed point free.

We claim that $g$ cannot invert each of $w_1, w_2$ and $w_3$ -- for if so,
then $g$ would invert each element in $G/N$ and $\langle w_1N, w_2N\rangle$
is isomorphic to $S_3$.  So choose a $w_i$ not inverted by $g$.
Then $g=(gw_i)w_i^{-1}$.  Since $gw_i$ does not have order $2$, it
is indecomposable and we have noted already that $w_i$ is indecomposable
and fixed point free.

This completes the proof.
\end{proof}

There are few irreducible groups containing transvections.
See \cite{Mc}.  If $G$ is a solvable primitive permutation group of degree $n$,
then $G$ does contain a regular normal subgroup.  Thus, using the previous
result and \cite{Mc} yields:

\begin{corollary} \label{solvable indecomposable}  If $G$ is
a primitive solvable subgroup of $S_n$, then one of the following
holds:
\begin{enumerate}
\item Every element of $G$ is indecomposable; 
\item $n=4$ and $G=S_4$; or 
\item $n=16$ and $G$ has a normal regular elementary
abelian subgroup $N$ of order $16$ and $G/N=\mathrm{0}_4^+(2)$.
 \end{enumerate}
\end{corollary}

We can now handle all primitive groups other than
the almost simple groups.

\begin{theorem} \label{productaction}
Assume that $G$ acts faithfully and primitively
on the set $X$ of cardinality $n > 1$. Assume
that $G$ is not almost simple. Every element of
$G$ can be written as a product of two indecomposable
elements, one of which is fixed point free.
\end{theorem}

\begin{proof} By the previous result, we may assume
that $G$ does not contain a regular normal subgroup.
We may also assume that some nontrivial element of
$G$ fixes at least $n/2$ points. It follows by the
structure of primitive groups \cite{AS}, the previous
result and \cite{GM} that $G$ preserves a Cartesian
product structure on $X$.

More precisely, we can write
$$
X=X_1 \times \ldots \times X_m,
$$
where $m > 1$, $|X_i|=e \ge 5$ and
$
G \le T:=S_e \wr S_m =W.S_m$ where
$$W = S_e \times \ldots \times S_e$$
acting coordinatewise
on $X$ and $S_m$ permutes the coordinates. Furthermore,
$G$ has a unique minimal normal subgroup
$$N:=L_1 \times \ldots \times L_m$$
where $L_i \cong L$ is a nonabelian simple and $L_i$ acts
on $X_i$ and trivially on $X_j$ for $j \ne i$.

Let $W_i$ be the $i$th copy of $S_e$ in $W$.

We claim that $g \in G$ is decomposable implies $g \in W_i$
for some $i$. It suffices to show that this is the case
for $T$. Suppose that $x,y \in T$ are nontrivial elements
and $M(x) \cap M(y)$ is empty. Suppose that $x$ acts on
an $X_i$ and $y$ on an $X_j$ with $j \ne i$. Choose $a \in X_i$
moved by $x$ and $b \in X_j$ moved by $y$. Then any point of
$X$ whose $i$th coordinate is $a$ and $j$th coordinate is
$b$ is moved by $x$ and $y$, a contradiction.

This shows that if $x$ and $y$ are both in $W$, then
they are both in $W_i$ for some $i$ and so also $xy$.
If neither $x$ nor $y$ is in $W$, then $x$ and $y$
each move at least $n-n/e > n/2$ points and so $M(x) \cap M(y)$
is nonempty. Finally, suppose that $x$ is not in $W$ and
$y \in W$. Arguing as above, we see that it suffices
to consider the case that $x$ permutes the $X_i$ transitively.
Say $y$ is nontrivial on $X_1$ and moves $a \in X_1$.
Then $x$ cannot fix all points of $X$ with first coordinate
$a$ and so $M(x) \cap M(y)$ is empty.

This proves the claim.

We now complete the proof of the result.

Let $g \in G$. If $g$ is not in $W$, then choose
$h \in N$ with $h$ not in $N \cap W_i=L_i$ for any
$i$ and $h$ fixed point free (just choose $h_1 \in L_1$
fixed point free and $h_2$ nontrivial). Then
$g=h(h^{-1}g)$ is the desired decomposition ($h^{-1}g$
is not in $W$ and so indecomposable). If $g \in W$,
we choose a similar $h$ guaranteeing that $h^{-1}g$
is not in $W_i$ for any $i$.

\end{proof}

\subsection{Almost Simple Groups}\label{almost-simple-groups}
We now consider almost simple groups. So
$G$ is an almost simple group and has socle $S$
and acts transitively on $X$ of cardinality $n > 1$.

We first deal with the cases $G=A_n$ or $S_n$.
Note that the lemma is just the theorem for
these groups.

\begin{lemma} \begin{enumerate}
\item Any element of $S_n$ can be written as a product
of an $n$-cycle and a $k$-cycle for some $k$.
\item If $n$ is even, then every element of $A_n$ can
be written as product $xy$ where $x$ has exactly
two orbits each of even length and $y$ is a $k$-cycle
or $y$ has precisely two nontrivial orbits each of even
length.
\end{enumerate}
\end{lemma}

\begin{proof}Suppose that $g$ has $k$ orbits.

Let $h$ be a $k$-cycle
moving precisely one point in each $g$-orbit. Then
$gh$ is an $n$-cycle, whence (1) holds.

Now suppose that $n$ is even and $g \in A_n$.
If $g=1$, the result is clear.
Otherwise,
write $g=xy$ where $x$ is an $n$-cycle and $y$
is a $k$-cycle. Necessarily $k$ is even and the construction
above shows that we can take $k < n$.

Let $t$ be a transposition moving at least $1$ point fixed by $y$.
Then $xt$ has precisely
$2$ orbits and we can pick $t$ so that each of the orbits
is even. Then $ty$ is either a $k+1$ cycle (if $t$ and $y$
are not disjoint) or has two nontrivial orbits (of length $2$ and $k$).
So $g=(xt)(ty)$, whence (2) holds.
\end{proof}

If no element fixes at least half the points, then clearly
every element is indecomposable.  By \cite{GM}, the only cases
to consider are dealt with in the next three lemmas.

\begin{lemma} \label{ksets}
Let $G=A_n$ or $S_n$ with $n \ge 5$
acting on $X$, the set of $k$-sets for some $k$ with $1 < k < n/2$.
Then every element of $G$ is indecomposable.
\end{lemma}

\begin{proof} We show that for $x,y$ nontrivial, $M(x)$ and
$M(y)$ have a nonempty intersection. Let $Y=\{1, 2, \ldots, n\}$.
If $x \in G$ and $j \in Y$, we write $xj$ for the image of $j$
under $x$.

First suppose that $x$ and $y$ move a common point in the natural
representation. So we may assume that $x$ and $y$ each move $1$.
Let $D$ be a $k$-set containing $1$ but missing $x1$ and $y1$.
Then $x$ and $y$ both move $D$.

Suppose that $x$ and $y$ move no common point in $Y$.
So we may assume that $x$ moves $1$ and $y$ moves
$2$. Let $D$ be a $k$-set containing $1,2$ but not containing
$x1$ and $y2$. Then $x$ and $y$ both move $D$.
\end{proof}

\begin{lemma} \label{sp} Let $G=\Sp(2d,2)$ with $d \ge 3$.
Let $X$ be the coset space $G/H$ where $H=\Or^{-}(2d,2)$
(note that this is the set of nondegenerate hyperplanes of
$-$ type in the $2d+1$ dimensional orthogonal module for $G$).
Every element of $G$ is indecomposable on $X$.
\end{lemma}

\begin{proof} Suppose that $M(x) \cap M(y)$ for $x,y$ nontrivial
in $G$. It is easy to see (cf \cite{GM}) that every nontrivial element
other than a transvection moves more than $|X|/2$ elements.
So we choose notation so that
$x$ is a transvection and $y \ne x$. Let $P=C_G(x)$.
Then $P$ is a maximal parabolic subgroup of $G$. Then $y$
fixes each coset of $H$ moved by $x$. The same is true
for any $P$-conjugate of $y$ and so $J:=\langle y^P \rangle$ does
as well. So $P$ normalizes $J$. Now $J$ is proper in $G$
and so as $G$ is simple and $P$ is maximal, $J$ is a nontrivial
normal subgroup of $P$. The subgroup generated by $x$ is
the unique minimal normal subgroup of $P$ and so $x \in J$.
However, $x$ certainly moves all the points of $M(x)$ and
this contradiction completes the proof.
\end{proof}

\begin{lemma} \label{or} Let $G^{\epsilon}=
\Or^{\epsilon}(2d,2)$ with $d > 2$.
Let $X$ be the set of singular vectors (if $\epsilon=-$)
or the set of nonsingular vectors (if $\epsilon=+)$.
Every element of $G^{\epsilon}$ is indecomposable on $X$.
\end{lemma}

\begin{proof} Let $J=\Sp(2d,2)$ and $Y$ the $J$-set described
in the previous lemma. Note that $G^{\epsilon}$ is a subgroup
of $J$ and so acts on $Y$. If $\epsilon=+$, then $Y \cong X$
at $G^+$-sets. Also, $G^{-}$ fixes one point of $Y$ and
the remaining orbit is isomorphic to $X$ as a $G^{-}$ set.
Thus, the result follows from the previous lemma.
\end{proof}

The previous three lemmas together with  \cite{GM} 
immediately yields:

\begin{theorem} Let $B$ be an almost simple group acting
primitively on $X$. Then either every element of $G$ is
indecomposable or $G$ contains $\alt(X)$.
\end{theorem}

For almost simple groups, we can weaken the assumption of primitivity.

\begin{theorem} Let $G$ be an almost simple group transitive permutation
group of degree $n$ and suppose that some element of $g$ is decomposable.
 Then $G$ is a symmetric group or alternating
group of degree $m$ for some $m$ dividing $n$.
\end{theorem}

\begin{proof} If $G$ is primitive on $X$, this
follows from the previous result. Suppose that $G$ is
not primitive on $X$ and some element $g \in G$ is decomposable
on $X$.  Write $g=g_1g_2$ where the $g_i$ are disjoint on $X$
(and each nontrivial).  Let $S$ be the socle of
$G$.  

We induct on $|X|$.
Let $Y=\{X_1, \ldots, X_t\}$ be a nontrivial $G$-invariant
partition of $X$ with $G$ primitive on $Y$.   Let
$K$ be the normal subgroup of $G$ acting trivially on $Y$.
If $K=1$, then $G$ is faithful and primitive  on $Y$, whence
$G=\alt(Y)$ or $\sym(Y)$.  Otherwise $S \le K$ (since it is
the unique minimal normal subgroup of subgroup of $G$ containing $S$).

Assume that $g_2$ is not in $K$.  
Choose notation so that $X_1, \ldots, X_s$ with $s > 1$ is an orbit
for $g_2$ and set $X'=X_1 \cup \ldots \cup X_s$.   Then
$g_2$ is fixed point free on this set and so $g_1$ must be trivial
on this set.  Since $S$ leaves $X'$-invariant, it follows that
the stabilizer of $X'$ acts faithfully on $X'$, a contradiction.

So we may assume that $g_1$ and $g_2$ are both trivial on $Y$,
whence they both act on $X_1$ and as above both act nontrivially
on $X_1$.  So by induction, the result follows.
\end{proof}

Note that the previous result actually gives more information
with a little more effort---when $G$ is an alternating or symmetric group, essentially
the only maximal subgroup containing $H$ is unique and is the
stabilizer of a point in the natural permutation representation
(being slightly careful when $m=6$).


Combining the results on almost simple groups allows us
to state a more precise version of Theorem \ref{productaction}.
Note that in the proof of that theorem, we saw that the only
decomposable elements were contained in a component $L$ of $G$
and in particular, the component would have to be a simple
group that admits an action with decomposable elements.
Indeed, it follows by \cite{AS} that this action corresponds
to a primitive action of $N_G(L)/C_G(L)$ and so by
the result on almost simple groups $L=A_d$.

Thus we have the following result that will be useful in the final section.

\begin{theorem} \label{indecomposable}  Let $G$ be a primitive
subgroup of $S_n$.  One of the following holds:
\begin{enumerate}
\item Every element of $G$ is indecomposable;
\item $G=A_n, n > 5$ or $S_n, n > 3$;
\item $n=d^t$ with $d \ge 5$ and $t \ge 2$, $G \le S_d \wr S_t$
and $G$ contains ${A_d}^{t}$;
\item $n=2^a, a > 2$, $G$ contains a regular normal elementary
abelian subgroup $N$ and $G=NH$ where $H$ is a point stabilizer
and $H$ is an irreducible subgroup of $\mathrm{Aut}(N)$ containing
transvections.
\end{enumerate}
\end{theorem}

\section{Permutation Polytopes}
Now let $G$ be any finite group, and let $\nu\colon
G\to\mathrm{GL}(\R^n)$ be a real representation. The {\em
representation polytope} associated with $\nu$ is the convex hull
of the image of $\nu$, a subset of
$\mathrm{End}_{\R}(\R^n)\approx\R^{n^2}$:
\[
P(\nu):=\mathrm{conv}\{\nu(g)\in\R^{n^2}\mid g\in G\}.
\]
For each $g\in G$, left multiplication by $\nu(g)$ defines a linear
automorphism of $\R^{n^2}$ sending $P(\nu)$ to itself and sending the
image of the identity element of $G$ to $\nu(g)$. Hence, the vertices
of $P(\nu)$ are precisely the images of elements of
$G$.

If $G$ a subgroup of the symmetric group,
$S_n$, we write $P(G)$ for $P(\nu_G)$ where $\nu_G$ is the
natural representation of $G$ as a group of $n\times n$
permutation matrices. In this case, we also identify each $g\in G$
with its image, $\nu(g)\in\R^{n^2}$. The polytope, $P(G)$, is called
the {\em permutation polytope} associated with the permutation group
$G$.

In this part of the paper, we establish two main results. First, we
show that as $G$ varies over subgroups of $S_n$, the corresponding
polytope has maximal dimension $(n-1)^2$ exactly when $G$ is
$2$-transitive. Next, we characterize some faces of $P(G)$ and
give a bound on the diameter of the edge graph of $P(G)$.

\subsection{Dimension}
We use the following standard theorem from representation theory:
\begin{theorem}[Frobenius and Schur \cite{CR},
\S27.8]\label{thm:frobenius-schur}
Let $G$ be a finite group, $K$ an algebraically closed field, and
$\rho_i\colon G\to\mathrm{GL}(K^{n_i})$ for $i=1,\dots,k$ a collection
of pairwise non-isomorphic irreducible matrix representations of $G$.
Let $x^{(r)}_{ij}$ denote the coordinate functions of $\rho_r$ for
each $r$. Then
the set $\{x^{(r)}_{ij}\}_{i,j,r}$ of all coordinate functions is
linearly independent over~$K$.
\end{theorem}

Let $\nu=\oplus\nu_i^{a_i}$ be the
irreducible decomposition of $\nu$ over the complex numbers.
\begin{theorem}\label{thm:representation-dimension}
The dimension of the representation polytope $P(\nu)$ is
\[
\dim P(\nu)=\sum_{\nu_i\neq 1}(\deg \nu_i)^2,
\]
the sum taken over all non-trivial components $\nu_i$, not counting
multiplicities.
\end{theorem}
\begin{proof} Let $\C[G]$ denote the group algebra, and let
$\nu_i$ be a representation of $G$ on a complex vector space $V_i$ for
each $i$.
There is a natural algebra homomorphism
\[
\Gamma_{\nu}\colon\C[G]\to\oplus_i\mathrm{End}_{\C}(V_i)^{a_i}\subset\mathrm{End}_{\C}(\C^n)
\]
determined by $g\mapsto \nu(g)$ for each $g\in G$ and extending
linearly. The mapping $\Gamma_{\nu}$ further factors through the
inclusion
\begin{eqnarray*}
\oplus_i\mathrm{End}_{\C}(V_i)&\to&\oplus_i\mathrm{End}_{\C}(V_i)^{a_i}\\
\oplus_i\phi_i&\mapsto&\oplus_i\phi_i^{a_i}
\end{eqnarray*}
where $\phi\in \mathrm{End}_{\C}(V_i)$ for each $i$.  The resulting
mapping of $\C[G]$ into $\oplus_{i=1}^k\mathrm{End}(V_i)$ is a
surjection by Theorem~\ref{thm:frobenius-schur}.

Restricting $\Gamma_{\nu}$ to $\R[G]$, the polytope $P(\nu)$ is the
convex hull of the image of $G$. Hence, the dimension of $P(\nu)$ will 
be the
dimension of the image of $\Gamma_\nu$ if the polytope contains the
zero vector in its affine span and will be one less, otherwise. So it
suffices to show that $P(\nu)$ does not contain $\vec{0}$ in its
affine span, $\mathrm{aff}(P(\nu))$, if and only if $\nu$ contains the
trivial representation as an irreducible factor. First, suppose
$\vec{0}\not\in\mathrm{aff}(P(\nu))$. The vector
$\textstyle\frac{1}{|G|}\sum_{g \in G}\nu(g)$ is an element of
$P(\nu)$, hence nonzero, and its linear span is clearly $G$-invariant;
thus, $\nu$ contains the trivial representation. Conversely, suppose
that $\nu$ contains the trivial representation. Then there exists a
nonzero $w\in\C^n$ such that $\nu(g)(w)=w$ for all $g\in G$. Given an
arbitrary element $x=\sum_{g\in G}a_g\nu(g)$ in
$\mathrm{aff}(P(\nu))$, we have $x(w)=(\sum a_g)w=w$, hence,
$x\neq\vec{0}$, as required.
\end{proof}

\begin{corollary}
If $\nu$ is a faithful representation, $P(\nu)$ is a simplex if and
only if each irreducible representation of $G$ appears up to
isomorphism as a component in the irreducible decomposition of $\nu$.
\end{corollary}
\begin{proof}
Let $\nu=\oplus \nu_i^{ a_i}$ be the irreducible decomposition of
$\nu$ over $\C$. The polytope $P(\nu)$ is a simplex if and only
if its dimension is one less then the number of vertices. In light of
Theorem~\ref{thm:representation-dimension}, the condition is equivalent
to
$|G|-1=\sum_{\nu_i\neq 1}(\deg \nu_i)^2$. However, a basic theorem of
representation theory says that $|G|=\sum_{\tau} (\dim \tau)^2$ where
the sum
is over a full set of representatives of the isomorphism classes of
irreducible
representations of $G$ (including the trivial representation).
\end{proof}
If $\nu$ is not faithful, let $H=\{g\in G\mid\nu(g)=1\}$. In this case,
$P(\nu)$ is a simplex if and only if the irreducible decomposition of
$\nu$ over $\C$ contains each irreducible representation of $G$
trivial on $H$.

\begin{corollary}\label{cor:transitivity} Let $G\leq S_n$ be a subgroup 
having $t$ orbits.
\begin{enumerate}
\item\label{cor:transitivity1} $\dim P(G)\leq (n-t)^2$ with equality
if and only if $\nu_G$ has
at most one non-trivial factor in its irreducible decomposition;
\item\label{cor:transitivity2} $\dim P(G)\leq (n-1)^2$ with equality
if and only if $G$ is $2$-transitive.
\item\label{cor:transitivity3} The dimension of the Birkhoff
polytope, $B_n$, is
$(n-1)^2$ for all $n\geq 1$.
\item\label{cor:transitivity4} The dimension of the polytope of even
permutation matrices, $A_n$, is $(n-1)^2$ for $n\geq 4$.
\end{enumerate}
\end{corollary}
\begin{proof}
Consider the irreducible decomposition of the permutation
representation $\nu_G=\oplus_i\nu_i^{a_i}$
over $\C$. It is well-known
from representation theory that the number of copies of the
trivial representation appearing in $\nu$ is the number of orbits,
$t$ (\cite{CR} \S 32.3). Let $\nu_1,\dots,\nu_k$ be the non-trivial
factors of $\nu_G$.
Then $\sum_{i=1}^k\deg\nu_i=n-t$ and by
Theorem~\ref{thm:representation-dimension}, the dimension of
$P(G)=\sum_{\nu_i\neq1}(\deg \nu_i)^2$. The sum is maximized when
$k\leq 1$. This proves part~\ref{cor:transitivity1}.

For part~\ref{cor:transitivity2}, by standard representation theory
of permutation groups, $G$ is $2$-transitive if and only if
$\nu_G=1+\tilde{\nu}_G$ for some irreducible $\tilde{\nu}_G$ (\cite{CR}
\S32.5).
Parts~\ref{cor:transitivity3} and~\ref{cor:transitivity4} then
follow since the relevant groups are $2$-transitive.
\end{proof}

\subsection{Faces}
Let $G\leq S_n$ be a permutation group, and identify elements of $G$
with $n\times n$ permutation matrices as usual.
For $g,h\in G$, write $h\preceq g$ if the set of cycles of $h$ is a
subset of the set of cycles of $g$ (so $M(h)\cap M(h^{-1}g)$ is empty).
The element $g$ is indecomposable when $h\preceq g$ always implies $h$
is
the identity or $g$.
\begin{theorem}\label{thm:smallest face} The smallest face of $P(G)$ 
containing $g,h\in G$ is
\[
F_{\{g,h\}}:=\mathrm{conv}\,\{hk\in G\mid k\preceq h^{-1}g\}.
\]
In particular, there is an edge connecting $g$ and $h$ if and only if
$h^{-1}g$ is indecomposable.
\end{theorem}
\begin{proof}
By symmetry, we may assume that $h$ is the identity, $e$, and show that
the smallest
face containing $g$ and $e$ is $\mathrm{conv}\{k\in G\mid k\preceq g\}$.
If $k\preceq g$, let $k'=k^{-1}g$. From
$g=kk'$ with $k,k'\preceq g$, it follows that
\begin{eqnarray}\label{pf:identity}
e+g=k+k'.
\end{eqnarray}
Let $c\in\R^{n^2}$ and $b\in\R$ with Euclidean inner products $\langle
c,g\rangle=\langle c,e\rangle=b$ and
$\langle c,f\rangle\leq b$ for all $f\in G$; so $c$ defines a face of
$P(G)$ containing $g$ and $e$. Equation~\ref{pf:identity} then
implies that $\langle c,k\rangle=\langle c,k'\rangle=b$, too. Hence,
any face containing $g$ and $e$ must also contain $k$ and $k'$.

For any matrix $m\in\R^{n^2}$, define the {\em support} of $m$ by
\[
\mathrm{supp}(m)=\{(i,j)\in \{1,\dots,n\}^2\mid m_{ij}\neq 0\}.
\]
Define the matrix $c\in\R^{n^2}$ by
\[
c_{ij}=\left\{
\begin{array}{cl}
1&\mbox{if $(i,j)\in\mathrm{supp}(g+e)$},\\
0&\mbox{otherwise.}
\end{array}\right.
\]
It follows that
$\langle c,g\rangle=\langle c,e\rangle=n$ and for any $f\in G$,
\[
\langle c,f\rangle=\sum_{(i,j)\in\mathrm{supp(g+e)}}f_{ij}\leq n
\]
with equality if and only if $f\preceq g$. Hence, $c$ defines a
face---the smallest face, $F_{\{g,e\}}$---containing both $g$ and $e$.
\end{proof}
Note that if $g=g_1\dots g_t$ with $g,g_1\cdots,g_t\in G$ and such that
the cycles of $g_1,\dots,g_t$ are disjoint, then
\[
g-e=\sum_{i=1}^t(g_i-e),
\]
hence, $g$ is affinely dependent on $g_1,\dots,g_t$.

A direct computation based on the theorem establishes the following
known results \cite{BR}, \cite{Yo}, \cite{BL}:
\begin{corollary}\
\begin{enumerate}
\item\label{cor:diam1} The diameter of $P(S_n)$ is $1$ for $n< 4$ and
is $2$ for $n\geq 4$.
\item\label{cor:diam2} The diameter of $P(A_n)$ is $1$ for $n<6$ and
is $2$ for $n\geq 6$.\
\end{enumerate}
\end{corollary}

Corollaries~\ref{cor:bound1} and~\ref{cor:bound2} translate into
bounds on the diameter of a permutation polytope.
\begin{corollary}\label{diameter} Let $G\leq S_n$.
The diameter of the polytope $P(G)$ is at
most $\min\{2t,\lfloor n/2\rfloor\}$, where $t$ is the number of 
nontrivial orbits of
$G$. In particular, if $G$ is transitive, the diameter of
$P(G)$ is at most $2$.
\end{corollary}
The bound is sharp. For example, take $G$ to be the direct product of
$t$
copies of the dihedral group on $4$ elements, naturally considered as
a subgroup of $S_{4t}$.

\section{Mixing Times}
In this section, we consider random walks on permutation polytopes
or equivalently on the Cayley graph of the permutation group $G$
with the corresponding generating set consisting of the
indecomposable elements of $G$. This problem was suggested
to us by Pak. The
question about the mixing time of random walks on
0-1 polytopes goes back some time.  See the survey
article \cite{Zi}.

 We generalize his result here.
First we recall some notation.
(see \cite{Pa}).

Let $G$ be a finite group and $S$ a symmetric generating set for
$G$ (i.e. $G = \langle S \rangle$ and $S=S^{-1}$). Let $Q^k(g)$
be the probability that a random product of $k$ elements of $S$ is
equal to $g$. Similarly, define $Q^k(A)$ to be the probability that
a random product of $k$ elements of $S$ is in the subset $A$ of $G$.
Let $U$ denote the uniform distribution on $G$.
Define the total variation distance,
$$
d(k):=(1/2)\sum_{g \in G} |Q^k(g) - 1/|G||= 
\max_{A \subseteq G} |Q^k(A) - U(A)|.
$$
So $d(k)$ measures how far the probability distribution $Q^k$ is
from the uniform distribution on $G$.

We now consider the case that $G$ is a subgroup of $S_n$ and $S$
is the set of indecomposable elements in $G$. Clearly, $S$ is
symmetric, $1 \in S$ and $G=\langle S \rangle$. We note that 
$Q^k \rightarrow U$ as $k \rightarrow \infty$ (i.e.
$d(k) \rightarrow 0$; this is standard since $S=S^{-1}$
and the Cayley graph is not bipartite -- see for example
\cite{AF}).  


\begin{theorem}\label{Pak1} Assume that $G$ is primitive of
degree $n$. If $G$ does not contain $A_n$, then $d(1) \rightarrow
0$ as $n \rightarrow \infty$.
In all cases, $d(2) \rightarrow 0$ as $n \rightarrow \infty$.
\end{theorem}

Pak \cite{Pa} proves this for the special case $G=S_n$.
The proof of this theorem follows easily from \S\ref{almost-simple-groups}
and Pak's result. Namely, by Theorem \ref{indecomposable}
one of the following holds:

\begin{enumerate}
\item $G=A_n$ or $S_n$;
\item $n=2^a$, $G$ contains a regular normal subgroup $N$ (elementary 
abelian
of order $2^a$) and a point stabilizer $H \le \mathrm{Aut}(N)$ contains 
transvections
and acts irreducibly on $N$; 
\item $n=d^t$ with $d \ge 5$, $t \ge 2$,
$G$ has a unique minimal normal subgroup $N=L \times \ldots \times L$
where $L \cong A_d$ and all decomposable elements of $G$ are contained
in one of the $t$ minimal normal subgroups of $N$; or
\item Every element of $G$ is indecomposable.
\end{enumerate}

First note, that if $d(1) \rightarrow 0$, it follows easily
that $d(2) \rightarrow 0$.  

In the first case, Pak \cite{Pa}
proved the result for $S_n$.  A trivial modification of his
proof shows that the result also holds for $A_n$.
As Pak points out, his proof used a well-known but unpublished
result of Lulov about the sum of the inverses of the degrees
of the irreducible representations of the symmetric groups.
A stronger version of this theorem is in Corollary 2.7 of
\cite{LS}.

Set $Y:=G \setminus{S}$.
So we only need prove that $|Y|/|G| \rightarrow
0$ as $n \rightarrow \infty$ in cases 2,3 and 4.

In the fourth case, $Y$ is empty.

Consider the second case.

In the  second case,
the only decomposable elements are fixed point free involutions (for 
they must be
the product of two elements each moving precisely $1/2$ the points and 
moving no
common points). Let $T$ be the set of involutions in $G$ which
have a fixed point and induce a transvection on $N$.  Note
that if $x \in T$, then $|xN \cap T|=2$ (indeed, 
$xN \cap T=x[x,N]$ and since $x$ acts as a transvection on
$N$, $|[x,N]|=2$).

The list of possible $H$ was determined by McLaughlin 
\cite{Mc}.
It follows easily from this that
$$
\lim_{a \rightarrow \infty} |T \cap H|/|G|^{1/2} = 0.
$$

Thus, $|Y| \le 4|T \cap H|^2$ and so $\lim_{a \rightarrow \infty}
|Y|/|G| \rightarrow 0$ as required.

Finally, consider the third case.  As we saw, the only decomposable
elements are in one of the $t$ normal subgroups of $N$.
Thus, $|Y| \le t(d!)$ and $|G| \ge (d!)^t$.  Since
$t > 1$, $|Y|/|G| \rightarrow 0$ as either $d$ or $t$ increases.

This completes the proof of the theorem.

We now give two examples to show that if $G$ is not primitive, the 
previous theorem
need not hold. More precisely, we produce a sequence of groups $G_p$ for 
$p$ an odd prime
such that for fixed $k$, $d(k)$ is bounded away from $0$. In the
first sequence,
the Cayley graph is close to bipartite and in the second sequence,
$Q^1$ is very small outside a proper normal subgroup.

Let $n=2p$. Let $x$ and $y$ be $p$-cycles in $S_n$ that are disjoint. Let
$u$ be an involution in $S_n$ with $uxu=y$. Set $G_p = \langle x, y, u 
\rangle$.
So $|G|=2p^2$ and has a normal elementary abelian subgroup $N:=\langle 
x, y \rangle$.
So $G$ is a transitive subgroup of $S_n$. Let $S$ be the set of 
indecomposable elements
in $G$.

Note that $xN \subset S$ and $N \cap S = \{ x^i, y^i | i = 0, 1, \ldots 
p-1\}$.
So $|S \cap N|=2p-1$. Thus, the probability that a random element of $S$ 
is in
$N$ is $(2p-1)/(p^2 + 2p-1) < 2/p$. In particular, we see that
$Q^k(N) > (1-2/p)^k$ if $k$ is even and $Q^k(xN) > (1-2/p)^k$ if $k$ is odd.
This shows that $d(k) \rightarrow 1/2$ as $p \rightarrow \infty$. 
In particular, 
the mixing time is unbounded.  Indeed, in the example,
we see that the mixing time is linear in $p$.

Pak \cite{Pa} did show that this could
happen for some $0,1$ polytopes---his example is 
essentially $\mathbb{Z}/2 \times S_n$.

We give another example that is similar in flavor to Pak's
example.  Let $J$ be a nonabelian group of order $qr$
with $q > r$ primes (so $r(q-1)$). Note that
$D$ embeds in $S_q$.  Let $p$ be a third
distinct prime and consider $G=\mathbb{Z}/p \wr J$
acting on $n:=pq$.  Let $N$ be the normal subgroup of
$G$ of index $r$.  Note that the number of indecomposable
 elements in $N$ is 
$(q-1)p^q + q(p-1) + 1$ while the number of indecomposable
elements outside $N$ is $(r-1)p^{q-1}$.  So the probability
that a random indecomposable element is not in  $N$ is less
than $1/p$.  Thus, the probability that a random product of
$k$ indecomposable elements is in $N$ is at least
$(1-1/p)^k$.  So for $p$ large compared to $k$, $Q^k$
is far from uniform.  Again, we see that the mixing time
is linear in $p$.

\end{document}